\documentclass[a4paper,12pt]{article}
\usepackage[francais]{babel}
\usepackage{amssymb, amsmath, eucal}
\usepackage{amsfonts}
\usepackage{amsthm}
\usepackage[all]{xy}
\usepackage{mathrsfs}
\usepackage[dvips]{graphics}
\usepackage[utf8]{inputenc}

\begin{document}

\newtheorem{The}{Theorem}[section]
\newtheorem{Lem}[The]{Lemma}
\newtheorem{Prop}[The]{Proposition}
\newtheorem{Cor}[The]{Corollary}
\newtheorem{Rem}[The]{Remark}
\newtheorem{Obs}[The]{Observation}
\newtheorem{SConj}[The]{Standard Conjecture}
\newtheorem{Titre}[The]{\!\!\!\! }
\newtheorem{Conj}[The]{Conjecture}
\newtheorem{Question}[The]{Question}
\newtheorem{Prob}[The]{Problem}
\newtheorem{Def}[The]{Definition}
\newtheorem{Not}[The]{Notation}
\newtheorem{Claim}[The]{Claim}
\newtheorem{Conc}[The]{Conclusion}
\newtheorem{Ex}[The]{Example}
\newtheorem{Fact}[The]{Fact}
\newcommand{\C}{\mathbb{C}}
\newcommand{\R}{\mathbb{R}}
\newcommand{\N}{\mathbb{N}}
\newcommand{\Z}{\mathbb{Z}}
\newcommand{\Q}{\mathbb{Q}}
\newcommand{\Proj}{\mathbb{P}}

\begin{center}

{\Large\bf The Albanese Map of sGG Manifolds and Self-Duality of the Iwasawa Manifold}

\end{center}

\begin{center}

{\large Dan Popovici}

\end{center}

\vspace{1ex}

\noindent {\small {\bf Abstract.} We prove that the three-dimensional Iwasawa manifold $X$, viewed as a locally holomorphically trivial fibration by elliptic curves over its two-dimensional Albanese torus, is self-dual in the sense that the base torus identifies canonically with its dual torus under a sesquilinear duality, the Jacobian torus of $X$, while the fibre identifies with itself. To this end, we derive elements of Hodge theory for arbitrary sGG manifolds, introduced in earlier joint work of the author with L. Ugarte, to construct in an explicit way the Albanese torus and map of any sGG manifold. These definitions coincide with the classical ones in the special K\"ahler and $\partial\bar\partial$ cases. The generalisation to the larger sGG class is made necessary by the Iwasawa manifold being an sGG, non-$\partial\bar\partial$, manifold. The main result of this paper can be seen as a complement from a different perspective to the author's very recent work where a non-K\"ahler mirror symmetry of the Iwasawa manifold was emphasised. We also hope that it will suggest yet another approach to non-K\"ahler mirror symmetry for different classes of manifolds.}

\vspace{1ex}

\section{Introduction}\label{section:introduction}

 In [Pop17], we proposed a new approach to the Mirror Symmetry Conjecture extended to {\bf possibly non-K\"ahler} compact complex manifolds. One of the main ideas was to substitute the Gauduchon cone for the classical K\"ahler cone that is empty on a non-K\"ahler manifold. The {\bf Iwasawa manifold}, a well-known compact non-K\"ahler manifold of complex dimension $3$ that was proved to have the weaker sGG property in [PU14] , was used in [Pop17] to illustrate our theory. The main result of [Pop17] was that the Iwasawa manifold is its own mirror dual. One of the arguments supporting this conclusion was the existence of a correspondence (that is holomorphic in the first argument, anti-holomorphic in the second) between a variation of Hodge structures (VHS) parametrised by what we called the {\it local universal family of essential deformations} of the Iwasawa manifold and a VHS parametrised by a subset of the {\it complexified Gauduchon cone} of this manifold.

\vspace{2ex}

 In the present paper, we give yet another criterion of a different nature by which the Iwasawa manifold is self-dual in a sesquilinear way. It states that in the well-known description of this manifold as a locally holomorphically trivial fibration by elliptic curves over a two-dimensional complex torus, both the base and the fibre are self-dual tori. This is the content of Theorem \ref{The:Iwasawa_mirror_self-duality} which is the main result of the paper. 

 The self-duality criterion is expressed in terms of the Albanese torus and map of the Iwasawa manifold that are manifestations of the Albanese torus and map (otherwise known to always be abstractly defined) we explicitly construct in full generality on any {\bf sGG manifold} by means of Hodge theory duly adapted to the specific context of possibly non-K\"ahler sGG manifolds. This construction occupies section \ref{section:Albanese-sGG}.

Our hope, motivating in part this note, is that the sesquilinear duality between the explicitly constructed Albanese torus and Jacobian torus of an arbitrary sGG manifold will show in the future how to guess the mirror dual of more general sGG manifolds that may not be mirror self-dual.  

\vspace{2ex} 

 Recall that the Iwasawa manifold $X = G/\Gamma$ is defined as the quotient of the Heisenberg group

$$G:=\left\{\begin{pmatrix}1 & z_1 & z_3\\   
                           0 & 1 & z_2\\
                   0 & 0 & 1\end{pmatrix}\,\, ; \,\, z_1, z_2, z_3\in\C\right\}\subset GL_3(\C)$$

\noindent  by its discrete subgroup $\Gamma\subset G$ of matrices with entries $z_1, z_2, z_3\in\Z[i]$. The map $(z_1,z_2,z_3)\mapsto (z_1,z_2)$ is easily seen to factor through the action of $\Gamma$ to define a locally holomorphically trivial proper holomorphic submersion

 \begin{equation}\label{introd_Alb-fibration}\pi : X\to B\end{equation}

\noindent whose base $B=\C^2/\Z[i]\oplus \Z[i]=\C/\Z[i]\times\C/\Z[i] $ is a two-dimensional complex torus (even an Abelian variety) and whose fibres are all isomorphic to the Gauss elliptic curve $\C/\Z[i]$. The torus $B$ and the map (\ref{introd_Alb-fibration}) are the Albanese torus, resp. Albanese map of the Iwasawa manifold in the standard sense in which these objects are associated with any compact complex manifold using a universal property (cf. e.g. [Uen75, chapter IV, $\S.9$]). 

 We give in section \ref{section:Albanese-sGG} a precise description of the Albanese torus and map that is valid on every sGG manifold (hence also on the Iwasawa manifold).

Recall that from the invariance under the action of $\Gamma$ of the $\C^3$-valued holomorphic $1$-form on $G$

$$G\ni M=\begin{pmatrix}1 & z_1 & z_3\\   
                           0 & 1 & z_2\\
                   0 & 0 & 1\end{pmatrix}  \mapsto M^{-1}\, dM = \begin{pmatrix}0 & dz_1 & dz_3-z_1\, dz_2\\   
                           0 & 0 & dz_2\\
                   0 & 0 & 0\end{pmatrix}$$ 
\noindent we get three holomorphic $1$-forms $\alpha,\beta, \gamma$ on the Iwasawa manifold induced respectively by the forms $dz_1, dz_2, dz_3-z_1dz_2$ on $\C^3$. They are such that

$$d\alpha = d\beta = 0 \hspace{2ex} \mbox{and}  \hspace{2ex} d\gamma = \partial\gamma = -\alpha\wedge\beta \neq 0 \hspace{3ex} \mbox{on}\hspace{1ex} X.$$

\noindent The forms $\alpha,\beta, \gamma$, that we call {\it structural}, and their conjugates are known to determine the whole cohomology of $X$ (cf. e.g. [Sch07]).

 Considering the Kuranishi family $(X_t)_{t\in\Delta}$ (that is known to be unobstructed by a result of Nakamura, although we shall not use this fact in the present paper) of the Iwasawa manifold $X=X_0$, it is known (cf. e.g. [Ang14, p. 75-77]) that there exist $C^{\infty}$ families $(\alpha_t)_{t\in\Delta}$, $(\beta_t)_{t\in\Delta}$, $(\gamma_t)_{t\in\Delta}$ of smooth $(1,\,0)$-forms on the fibres $(X_t)_{t\in\Delta}$ such that $\alpha_0=\alpha$, $\beta_0=\beta$ and $\gamma_0=\gamma$ and such that the forms $\alpha_t,\beta_t, \gamma_t$ and their conjugates determine the whole cohomology of $X_t$ (cf. e.g. [Ang14, p. 77-84]). 

\vspace{2ex}

 We will exploit the fact that the structural forms $\alpha_t, \beta_t, \gamma_t$, their conjugates and appropriate products thereof define {\bf canonical} bases in all the cohomology groups that we are interested in on every $X_t$ with $t$ sufficiently close to $0$. This will allow us to deduce from the general explicit construction in section \ref{section:Albanese-sGG} that the Albanese torus $\mbox{Alb}(X_t)$ of any small deformation $X_t$ of the Iwasawa manifold $X_0$ is {\bf self-dual} (cf. Lemma \ref{Lem:identification_tori}). Theorem  \ref{The:Iwasawa_mirror_self-duality} follows easily from this.

\section{The Albanese torus and map of an sGG manifold}\label{section:Albanese-sGG}

 Let $X$ be a compact complex manifold with $\mbox{dim}_{\C}X=n$.

\subsection{Elements of Hodge theory of $\partial\bar\partial$-manifolds} Recall that $X$ is said to be a {\bf $\partial\bar\partial$-manifold} if the $\partial\bar\partial$-lemma holds on $X$. This means that for every $p,q=0,1, \dots , n$ and for every $d$-closed smooth $(p,\,q)$-form $u$ on $X$, the following exactness conditions are equivalent\!:

\begin{equation}\label{eqn:ddbar-def}u\in\mbox{Im}\,d \iff u\in\mbox{Im}\,\partial \iff u\in\mbox{Im}\,\bar\partial \iff u\in\mbox{Im}\,\partial\bar\partial.\end{equation}

\noindent It is well known (see e.g. [Pop14] for a rundown on the basic properties of these manifolds) that on any $\partial\bar\partial$-manifold, the Hodge decomposition and the Hodge symmetry hold in the following sense\!\!: there exist {\bf canonical} (i.e. depending only on the complex structure of $X$) isomorphisms

\begin{equation}\label{eqn:Hodge-decomp}H^k_{DR}(X,\,\C)\simeq\bigoplus\limits_{p+q=k}H^{p,\,q}_{\bar\partial}(X,\,\C) \hspace{2ex} \mbox{and} \hspace{2ex} H^{p,\,q}_{\bar\partial}(X,\,\C)\simeq\overline{H^{q,\,p}_{\bar\partial}(X,\,\C)}, \hspace{2ex} k=0,1, \dots , 2n,\end{equation}

\noindent where $H^k_{DR}(X,\,\C)$ stands for the De Rham cohomology group of degree $k$, while $H^{p,\,q}_{\bar\partial}(X,\,\C)$ stands for the Dolbeault cohomology group of bidegree $(p,\,q)$. The inverse of the former isomorphism and the latter isomorphism are respectively defined by

$$([u^{p,\,q}]_{\bar\partial})_{p+q=k}\mapsto \bigg\{\sum\limits_{p=q=k}u^{p,\,q}\bigg\}_{DR}, \hspace{2ex}  [u]_{\bar\partial}\mapsto \overline{[\bar{u}]_{\bar\partial}}.$$ 

\noindent This is made possible by the fact that the $\partial\bar\partial$-lemma ensures the existence of a $d$-closed representative in {\it every} Dolbeault cohomology class $[u]_{\bar\partial}$ of any bidegree $(p,\,q)$ (see e.g. [Pop13, Lemma 3.1]). It also ensures that the above maps are independent of the choice of $d$-closed representatives in the classes involved. The $\partial\bar\partial$-lemma also defines {\it canonical} isomorphisms between any two of the cohomology groups $H^{p,\,q}_{BC}(X,\,\C)$ (Bott-Chern), $H^{p,\,q}_{\bar\partial}(X,\,\C)$ (Dolbeault) and $H^{p,\,q}_A(X,\,\C)$ (Aeppli), so in particular the Hodge decomposition (\ref{eqn:Hodge-decomp}) holds with any of $H^{p,\,q}_{BC}(X,\,\C)$ and $H^{p,\,q}_A(X,\,\C)$ in place of $H^{p,\,q}_{\bar\partial}(X,\,\C)$. 

 In other words, $\partial\bar\partial$-manifolds behave cohomologically as compact K\"ahler manifolds do. In particular, the {\bf Jacobian} and {\bf Albanese tori} and {\bf maps} can be defined on $\partial\bar\partial$-manifolds in a way identical to the one they are defined on compact K\"ahler manifolds.

\subsection{Elements of Hodge theory of sGG manifolds}

 The first purpose of this paper is to show that the {\bf Jacobian} and {\bf Albanese tori} and {\bf maps} can still be defined using Hodge theory in the larger class of {\bf sGG manifolds} (cf. [PU14]) with only minor modifications of the construction from the $\partial\bar\partial$ case. We will show that this is possible despite the fact that sGG manifolds need not admit a Hodge decomposition with symmetry in the standard sense of (\ref{eqn:Hodge-decomp}), but only a much weaker version thereof (cf. the splittings (\ref{eqn:H1-decomp-sGG}) and (\ref{eqn:H2n-1-decomp-sGG}) below that will play a key role in the sequel and what was called a {\it fake Hodge decomposition} in [PU14] that will not be used in this paper).

The {\bf sGG class} of compact complex manifolds, introduced in [PU14], strictly contains the class of $\partial\bar\partial$-manifolds, the best known example of an sGG manifold that is not a $\partial\bar\partial$-manifold being the {\bf Iwasawa manifold}. Recall the following equivalences (cf. [PU14])\!\!:

\begin{eqnarray}\label{eqn:sGG-characterisations}\nonumber X \hspace{1ex} \mbox{is sGG} & \stackrel{(a)}{\iff} & {\cal SG}_X = {\cal G}_X  \stackrel{(b)}{\iff}  \mbox{every Gauduchon metric on}\,\,X \,\, \mbox{is strongly Gauduchon}\\
\nonumber & \stackrel{(c)}{\iff} & \forall u\in C^{\infty}_{n,\,n-1}(X,\,\C)\cap\ker d, \hspace{1ex} \mbox{the implication holds\!\!:} \hspace{1ex} u\in\mbox{Im}\,\partial \implies u\in\mbox{Im}\,\bar\partial\\
\nonumber & \stackrel{(d)}{\iff} & b_1 = 2\,h^{0,\,1}_{\bar\partial},\end{eqnarray}

\noindent where $(a)$ is the definition (given in [PU14]) of sGG manifolds requiring the sG cone ${\cal SG}_X$ of $X$ to equal the (a priori larger) Gauduchon cone ${\cal G}_X$ (see [Pop15] for the terminology), $(b)$ is easily seen to be equivalent to $(a)$ (see e.g. [Pop14] for a reminder of the terminology), $(c)$ expresses the sGG property as a special case of the $\partial\bar\partial$-lemma (cf. [Pop15, Observation 5.3] --- the reader unfamiliar with the terminology of the other equivalences may wish to take equivalence $(c)$ as the definition of sGG manifolds), while $(d)$ is one of the numerical characterisations proved in [PU14]. Actually $b_1 \leq 2\,h^{0,\,1}_{\bar\partial}$ on every compact complex manifold and the equality characterises the sGG manifolds ([PU14, Theorem 1.5]). 

Moreover, by [PU14, Theorem 3.1], on every compact complex manifold $X$, the following canonical linear map\!\!:

\begin{eqnarray}\label{eqn:H1-decomp-sGG} F\,:\,H^1_{DR}(X,\,\C) \longrightarrow  H^{0,\,1}_{\bar\partial}(X,\,\C)\oplus\overline{H^{0,\,1}_{\bar\partial}(X,\,\C)},\hspace{3ex} F(\{\alpha\}_{DR}) := ([\alpha^{0,\,1}]_{\bar\partial},\,\overline{[\overline{\alpha^{1,\,0}}]_{\bar\partial}}),\end{eqnarray}

\noindent is well defined and injective. Furthermore, $X$ is sGG if and only if $F$ is an isomorphism. Equivalently, the dual linear map

\begin{eqnarray}\label{eqn:H2n-1-decomp-sGG} F^{\star}:\,H^{n,\,n-1}_{\bar\partial}(X,\,\C)\oplus\overline{H^{n,\,n-1}_{\bar\partial}(X,\,\C)} \longrightarrow  H^{2n-1}_{DR}(X,\,\C),\hspace{3ex} F^{\star}([\beta]_{\bar\partial},\,\overline{[\gamma]}_{\bar\partial}) := \{\beta + \bar\gamma\}_{DR},\end{eqnarray}

\noindent is surjective for any $X$, while $X$ is sGG if and only if $F^{\star}$ is an isomorphism.

 Thus, the canonical splittings (\ref{eqn:H1-decomp-sGG}) and (\ref{eqn:H2n-1-decomp-sGG}) of $H^1_{DR}(X,\,\C)$ and resp. $H^{2n-1}_{DR}(X,\,\C)$ are the weaker substitutes for the Hodge decomposition (\ref{eqn:Hodge-decomp}) in degrees $1$, resp. $2n-1$, afforded to sGG manifolds. Clearly, when $X$ is a $\partial\bar\partial$-manifold, (\ref{eqn:H1-decomp-sGG}) and (\ref{eqn:H2n-1-decomp-sGG}) coincide with the splittings for $k=1$, resp. $k=2n-1$, in (\ref{eqn:Hodge-decomp}).

\begin{Cor}\label{Cor:H01_injection_H1} For every {\bf sGG manifold} $X$, the Dolbeault cohomology group $H^{0,\,1}_{\bar\partial}(X,\,\C)$ {\bf injects canonically} into the De Rham cohomology group $H^1_{DR}(X,\,\C)$. The canonical injection $j:H^{0,\,1}_{\bar\partial}(X,\,\C)\hookrightarrow H^1_{DR}(X,\,\C)$ is obtained as the composition of the injective linear maps

$$H^{0,\,1}_{\bar\partial}(X,\,\C)\hookrightarrow H^{0,\,1}_{\bar\partial}(X,\,\C)\oplus\overline{H^{0,\,1}_{\bar\partial}(X,\,\C)} \stackrel{F^{-1}}{\longrightarrow} H^1_{DR}(X,\,\C).$$

\end{Cor}

\noindent {\it Proof.} The sGG assumption ensures that the canonical linear map $F$ defined in (\ref{eqn:H1-decomp-sGG}) is an isomorphism. Then so is its inverse $F^{-1}$.   \hfill $\Box$

\vspace{3ex}

 The canonical splittings (\ref{eqn:H1-decomp-sGG}) and (\ref{eqn:H2n-1-decomp-sGG}) enable one to construct canonically and explicitly the {\it Jacobian variety} (cf. Definition \ref{Def:jacobian-variety}) and the {\it Albanese variety} (cf. Definition \ref{Def:albanese-variety}) of any sGG manifold by imitating the classical constructions on compact K\"ahler manifolds with the necessary modifications. The details are spelt out in $\S.$\ref{subsection:jacobian} and  $\S.$\ref{subsection:albanese}.

\subsection{The Jacobian variety of an sGG manifold}\label{subsection:jacobian}

 Let $X$ be an sGG manifold with $\mbox{dim}_{\C}X=n$. The inclusions $\Z\subset\R\subset\C\subset{\cal O}$ induce morphisms 

$$H^1(X,\,\Z)\longrightarrow H^1(X,\,\R)\longrightarrow H^1(X,\,\C)\longrightarrow H^1(X,\,{\cal O})$$

\noindent where the image of $H^1(X,\,\Z)$ is a lattice in $H^1(X,\,\R)$. On the other hand, the map $H^1(X,\,\R)\rightarrow H^{0,\,1}_{\bar\partial}(X,\,\C)$ obtained by composing the maps $H^1(X,\,\R)\rightarrow H^1(X,\,\C)\rightarrow H^1(X,\,{\cal O})\simeq H^{0,\,1}_{\bar\partial}(X,\,\C)$ identifies canonically with the composite map 

$$H^1_{DR}(X,\,\R)\stackrel{j_1}{\hookrightarrow} H^1_{DR}(X,\,\C)\stackrel{p_1\circ F}\longrightarrow H^{0,\,1}_{\bar\partial}(X,\,\C),$$

\noindent where $j_1$ is the natural injection and $p_1\,\,:\,\,H^{0,\,1}_{\bar\partial}(X,\,\C)\oplus\overline{H^{0,\,1}_{\bar\partial}(X,\,\C)}\longrightarrow H^{0,\,1}_{\bar\partial}(X,\,\C)$ is the projection onto the first factor. Since $F$ is an {\bf isomorphism} (thanks to $X$ being {\bf sGG}), we get that 

$$p_1\circ F\circ j_1\,\,:\,\,H^1_{DR}(X,\,\R)\longrightarrow H^{0,\,1}_{\bar\partial}(X,\,\C)$$

\noindent is an isomorphism. Hence $\mbox{Im}\,H^1(X,\,\Z)$ is a lattice in $H^{0,\,1}_{\bar\partial}(X,\,\C)$. As a result, we can put

\begin{Def}\label{Def:jacobian-variety} The {\bf Jacobian variety} of an $n$-dimensional sGG manifold $X$ is defined exactly as in the K\"ahler case as the $q$-dimensional complex torus

\begin{equation}\label{eqn:Jac-def}\mbox{Jac}(X):=H^{0,\,1}_{\bar\partial}(X,\,\C)/\mbox{Im}\,H^1(X,\,\Z),\end{equation}

\noindent where $q:=h^{0,\,1}_{\bar\partial}(X)$ stands for the irregularity of $X$.

\end{Def}

\subsection{The Albanese variety of an sGG manifold} \label{subsection:albanese}

  Let once again $X$ be an sGG manifold with $\mbox{dim}_{\C}X=n$. In a way similar to the above discussion, we have morphisms

$$H^{2n-1}(X,\,\Z)\longrightarrow H^{2n-1}(X,\,\R)\stackrel{j_{2n-1}}{\longrightarrow} H^{2n-1}(X,\,\C)\stackrel{(F^{\star})^{-1}}{\longrightarrow} H^{n,\,n-1}_{\bar\partial}(X,\,\C)\oplus\overline{H^{n,\,n-1}_{\bar\partial}(X,\,\C)},$$

\noindent where $\mbox{Im}\,H^{2n-1}(X,\,\Z)$ is a lattice in $H^{2n-1}(X,\,\R)$ (a general feature of any compact complex manifold $X$) and $(F^{\star})^{-1}$ is an {\bf isomorphism} (thanks to $X$ being {\bf sGG}). If we denote by $p_2\,\,:\,\,H^{n,\,n-1}_{\bar\partial}(X,\,\C)\oplus\overline{H^{n,\,n-1}_{\bar\partial}(X,\,\C)}\longrightarrow \overline{H^{n,\,n-1}_{\bar\partial}(X,\,\C)}$ the projection onto the second factor, then 

$$p_2\circ (F^{\star})^{-1}\circ j_{2n-1}\,\,:\,\,H^{2n-1}_{DR}(X,\,\R)\longrightarrow \overline{H^{n,\,n-1}_{\bar\partial}(X,\,\C)}$$ 

\noindent is an isomorphism and therefore $\mbox{Im}\,H^{2n-1}(X,\,\Z)$ is a lattice in $\overline{H^{n,\,n-1}_{\bar\partial}(X,\,\C)}\simeq (\overline{H^{0,\,1}_{\bar\partial}(X,\,\C)})^{\star}$.

\begin{Def}\label{Def:albanese-variety} The {\bf Albanese variety} of an $n$-dimensional sGG manifold $X$ is the complex torus

\begin{equation}\label{eqn:Alb-def}\mbox{Alb}(X):= \overline{H^{n,\,n-1}_{\bar\partial}(X,\,\C)}/\mbox{Im}\,H^{2n-1}(X,\,\Z) = \bigg(\overline{H^{0,\,1}_{\bar\partial}(X,\,\C)}\bigg)^{\star}/\mbox{Im}\,H^1(X,\,\Z)^\star.\end{equation}

\end{Def}

The spaces $H^{n,\,n-1}_{\bar\partial}(X,\,\C)$ and $H^{0,\,1}_{\bar\partial}(X,\,\C)$ are dual under the Serre duality, while $H^{2n-1}(X,\,\Z)$ and $H^1(X,\,\Z)$ are Poincar\'e dual.

 Recall that in the standard case when $X$ is K\"ahler, the Albanese torus of $X$ is defined as the quotient

$$H^{n-1,\,n}(X,\,\C)/\mbox{Im}\,H^{2n-1}(X,\,\Z).$$ 

\noindent Since, by Hodge symmetry, the conjugation defines an isomorphism $H^{n-1,\,n}_{\bar\partial}(X,\,\C)\simeq \overline{H^{n,\,n-1}_{\bar\partial}(X,\,\C)}$ when $X$ is K\"ahler, our Definition \ref{Def:albanese-variety} of the Albanese torus coincides with the standard defintion in the K\"ahler case.

\vspace{3ex}

\begin{Conc}\label{Conc:dual_Jacobi-Albanese} We can now conclude from Definitions \ref{Def:jacobian-variety} and \ref{Def:albanese-variety} that the {\bf Jacobian torus} and the {\bf Albanese torus} of any sGG manifold $X$ are {\bf dual tori} in the sense of the following {\bf sesquilinear duality} obtained by composing the bilinear Serre duality with the conjugation in the second factor\!\!:

\begin{equation}\label{eqn:sesquilinear_Serre}H^{0,\,1}_{\bar\partial}(X,\,\C)\times\overline{H^{n,\,n-1}_{\bar\partial}(X,\,\C)}\longrightarrow\C, \hspace{3ex} ([\alpha]_{\bar\partial},\,\overline{[\beta]}_{\bar\partial})\mapsto\int\limits_X\alpha\wedge\beta.\end{equation}

\end{Conc}

\subsection{The Albanese map of an sGG manifold}\label{subsection:Albanese-map_sGG}

 We can now easily adapt to the general context of sGG manifolds $X$ the construction of the Albanese map $\alpha:X\longrightarrow\mbox{Alb}(X)$ from the familiar K\"ahler case. We shall follow the presentation and use the notation of [Dem97, $\S.9.2$].

 Let $X$ be an sGG manifold with $\mbox{dim}_{\C}X=n$. The standard isomorphism 

$$H_1(X,\,\Z)\longrightarrow H^{2n-1}(X,\,\Z)$$

\noindent given by the Poincar\'e duality is induced by the map $[\xi]\mapsto\{I_\xi\}_{DR}\in H^{2n-1}_{DR}(X,\,\R)$ associating with the homology class $[\xi]$ of every loop $\xi$ in $X$ the De Rham cohomology class of the current of integration $I_\xi$ over $\xi$. Using this isomorphism, the expression (\ref{eqn:Alb-def}) of the Albanese torus of $X$ transforms to

\begin{equation}\label{eqn:Alb_1}\mbox{Alb}(X)= \bigg(\overline{H^{0,\,1}_{\bar\partial}(X,\,\C)}\bigg)^{\star}/\mbox{Im}\,H_1(X,\,\Z),\end{equation}

\noindent where the map $H_1(X,\,\Z)\longrightarrow \overline{H^{0,\,1}_{\bar\partial}(X,\,\C)}^{\star}$ is defined by 

\begin{equation}\label{eqn:I-tilde_xi}[\xi]\mapsto \widetilde{I}_\xi:=\bigg(\overline{[v]}\mapsto\int\limits_\xi \overline{\{v\}}\bigg), \hspace{3ex} \mbox{where} \hspace{1ex} \{v\}:=j([v])\in H^1_{DR}(X,\,\C).\end{equation}

\noindent We have used the canonical injection $j:H^{0,\,1}_{\bar\partial}(X,\,\C)\hookrightarrow H^1_{DR}(X,\,\C)$ defined in Corollary \ref{Cor:H01_injection_H1} and the fact that the integral $\int_\xi \overline{\{v\}}$ depends only on the homology class $[\xi]$ and on the cohomology class $\overline{\{v\}}$ (so not on the actual representatives of these classes).

\begin{Def}\label{Def:albanese-map_def} Let $X$ be an sGG manifold. Fix a base point $a\in X$. For every point $x\in X$, let $\xi$ be any path from $a$ to $x$ and let $\widetilde{I}_\xi\in\overline{H^{0,\,1}_{\bar\partial}(X,\,\C)}^{\star}$ be the linear functional defined in (\ref{eqn:I-tilde_xi}).   

 The canonical holomorphic map

\begin{equation}\label{eqn:albanese-map_def}\nonumber\alpha: X\longrightarrow\mbox{Alb}(X)= \bigg(\overline{H^{0,\,1}_{\bar\partial}(X,\,\C)}\bigg)^{\star}/\mbox{Im}\,H_1(X,\,\Z), \hspace{3ex} x\mapsto \widetilde{I}_\xi \hspace{2ex} \mbox{mod} \hspace{1ex} \mbox{Im}\,H_1(X,\,\Z),\end{equation}

\noindent will be called the {\bf Albanese map} of the sGG manifold $X$.

\end{Def}

 Note that the class of $\widetilde{I}_\xi$ modulo $\mbox{Im}\,H_1(X,\,\Z)$ does not depend on the choice of path $\xi$ from $a$ to $x$ because for any other such path $\eta$, $\widetilde{I}_{\eta^{-1}\,\xi}\in\mbox{Im}\,H_1(X,\,\Z)$. Also note that definition (\ref{eqn:albanese-map_def}) of the Albanese map for sGG manifolds $X$ coincides with the standard definition when $X$ is K\"ahler. Indeed, in the K\"ahler case, $\overline{H^{0,\,1}_{\bar\partial}(X,\,\C)}$ is canonically isomorphic to $H^{1,\,0}_{\bar\partial}(X,\,\C)$ by Hodge symmetry. Moreover, the role played by the canonical injection $j:H^{0,\,1}_{\bar\partial}(X,\,\C)\hookrightarrow H^1_{DR}(X,\,\C)$ defined in Corollary \ref{Cor:H01_injection_H1} when $X$ is sGG is an apt substitute for the fact that every holomorphic $1$-form (i.e. the unique representative of every element in $H^{1,\,0}_{\bar\partial}(X,\,\C)$) is $d$-closed when $X$ is K\"ahler or merely $\partial\bar\partial$.

 As in the standard K\"ahler case, we have an alternative description of the Albanese map.

\begin{Obs}\label{Obs:alternative-description_Albanese-map} Let $X$ be an sGG manifold with $\mbox{dim}_\C X=n$. Using the expression (\ref{eqn:Alb-def}) of the Albanese torus of $X$, the Albanese map of $X$ is given by

\begin{equation}\label{eqn:albanese-map_alternative}\nonumber\alpha: X\longrightarrow\mbox{Alb}(X)= \overline{H^{n,\,n-1}_{\bar\partial}(X,\,\C)}/\mbox{Im}\,H^{2n-1}(X,\,\Z), \hspace{3ex} x\mapsto \overline{\{I_\xi\}^{n,\,n-1}} \hspace{2ex} \mbox{mod} \hspace{1ex} \mbox{Im}\,H^{2n-1}(X,\,\Z),\end{equation}

\noindent where $ \overline{\{I_\xi\}^{n,\,n-1}}\in\overline{H^{n,\,n-1}_{\bar\partial}(X,\,\C)}$ is the projection of the De Rham cohomology class $\{I_\xi\}_{DR}\in H^{2n-1}_{DR}(X,\,\R)$ onto $\overline{H^{n,\,n-1}_{\bar\partial}(X,\,\C)}$ w.r.t. the isomorphism

$$(F^{\star})^{-1}:\,H^{2n-1}_{DR}(X,\,\C)\stackrel{\simeq}{\longrightarrow} H^{n,\,n-1}_{\bar\partial}(X,\,\C)\oplus\overline{H^{n,\,n-1}_{\bar\partial}(X,\,\C)}  $$

\noindent induced by (\ref{eqn:H2n-1-decomp-sGG}). As usual, $I_\xi$ stands for the current of integration over the path $\xi$ from $a$ to $x$ in $X$.

\end{Obs}

 Note that in Observation \ref{Obs:alternative-description_Albanese-map} the only difference in the sGG case compared with the standard K\"ahler (or $\partial\bar\partial$) case is the substitution of $\overline{H^{n,\,n-1}_{\bar\partial}(X,\,\C)}$ for $H^{n-1,\,n}_{\bar\partial}(X,\,\C)$. These spaces are isomorphic by Hodge symmetry when $X$ is K\"ahler or merely $\partial\bar\partial$.

\section{Application to the mirror self-duality of the sGG Iwasawa manifold}\label{section:application_mirror-Iwasawa}

 In this section, we apply the construction in $\S.$\ref{section:Albanese-sGG} to the Iwasawa manifold that is known to not be a $\partial\bar\partial$-manifold (see e.g. [Pop14]). However, the Iwasawa manifold $X=X_0$ and all its small deformations in its Kuranishi family $(X_t)_{t\in\Delta}$ are sGG compact complex manifolds of dimension $3$ (cf. [PU14]). So, the extension to the sGG context of the classical constructions of the Albanese torus and map from the $\partial\bar\partial$ case, performed in $\S.$\ref{subsection:jacobian} and $\S.$\ref{subsection:albanese}, is key to our purposes here.

 For the Iwasawa manifold $X=X_0$ and all its small deformations $(X_t)_{t\in\Delta}$, the Albanese maps

$$\pi_t:X_t\longrightarrow \mbox{Alb}(X_t):=B_t,  \hspace{3ex} t\in\Delta,$$

\noindent have simple explicit descriptions and $\pi:=\pi_0:X_0\to B_0$ is a locally holomorphically trivial fibration whose fibre $\pi^{-1}(s)$ is the Gauss elliptic curve $\C/\Z[i]$ and whose base is the $2$-dimensional complex torus $\C/\Z[i]\times\C/\Z[i]$.

\vspace{3ex}

 First, we show that the Albanese torus of every small deformation $X_t$ of the Iwasawa manifold $X=X_0$ is {\bf self-dual} in the context of the construction of section \ref{section:Albanese-sGG}.

\begin{Lem}\label{Lem:identification_tori} Let $(X_t)_{t\in\Delta}$ be the Kuranishi family of the Iwasawa manifold $X=X_0$. Thus $n=\mbox{dim}_\C X_t = 3$. For every $t\in\Delta$ sufficiently close to $0$, the dual Jacobian and Albanese tori $\mbox{Jac}(X_t)$ and $\mbox{Alb}(X_t)$ can be identified {\bf canonically} in the following sense.

 There exist {\bf canonical} isomorphisms

\begin{equation}\label{eqn:identification_tori}H^{0,\,1}_{\bar\partial}(X_t,\,\C)\simeq H^{3,\,2}_{\bar\partial}(X_t,\,\C) \hspace{2ex} \mbox{and} \hspace{2ex} H^1(X_t,\,\Z)\simeq H^5(X_t,\,\Z), \hspace{3ex} t\in\Delta.\end{equation}

\end{Lem}

\noindent {\it Proof.} Dual finite-dimensional vector spaces are, of course, isomorphic, so the main feature of the isomorphisms (\ref{eqn:identification_tori}) is their canonical nature. By ``canonical'' we mean ``depending only on the complex or differential structure, independent of any choice of metric''. As can be seen below, the canonical nature of these isomorphisms follows from the existence of canonical bases, defined by the structural differential forms $\alpha_t,\beta_t, \gamma_t$ mentioned in the introduction and their conjugates, in the vector spaces involved.

 From [Sch07, p.6] and [Ang14, $\S.2.2.2$, $\S.2.2.3$], we gather that the vector spaces featuring in (\ref{eqn:identification_tori}) are generated by the structural $(1,\,0)$-forms $\alpha_t, \beta_t, \gamma_t$ as follows:

\begin{eqnarray}\label{eqn:generation}\nonumber H^{0,\,1}_{\bar\partial}(X_t,\,\C) & = & \bigg\langle[\bar\alpha_t]_{\bar\partial},\,[\bar\beta_t]_{\bar\partial}\bigg\rangle, \hspace{3ex} H^{3,\,2}_{\bar\partial}(X_t,\,\C) = \bigg\langle[\alpha_t\wedge\beta_t\wedge\gamma_t\wedge\bar\alpha_t\wedge\bar\gamma_t]_{\bar\partial},\,[\alpha_t\wedge\beta_t\wedge\gamma_t\wedge\bar\beta_t\wedge\bar\gamma_t]_{\bar\partial}\bigg\rangle,\\
H^1_{DR}(X_t,\,\C) & = & \bigg\langle\{\alpha_t\},\, \{\beta_t\},\, \{\bar\alpha_t\},\, \{\bar\beta_t\}\bigg\rangle, \end{eqnarray}
\begin{eqnarray}\nonumber  H^5_{DR}(X_t,\,\C) & = & \bigg\langle\{\alpha_t\wedge\beta_t\wedge\gamma_t\wedge\bar\alpha_t\wedge\bar\gamma_t\},\, \{\alpha_t\wedge\beta_t\wedge\gamma_t\wedge\bar\beta_t\wedge\bar\gamma_t\},\, \{\alpha_t\wedge\gamma_t\wedge\bar\alpha_t\wedge\bar\beta_t\wedge\bar\gamma_t\},\, \\
\nonumber & & \hspace{58ex} \{\beta_t\wedge\gamma_t\wedge\bar\alpha_t\wedge\bar\beta_t\wedge\bar\gamma_t\}\bigg\rangle,\end{eqnarray}

\noindent where $\{\,\,\,\}$ stands for De Rham cohomology classes. 

 Thus, the isomorphism $H^{0,\,1}_{\bar\partial}(X_t,\,\C)\simeq H^{3,\,2}_{\bar\partial}(X_t,\,\C)$ of (\ref{eqn:identification_tori}) is canonically defined by $[\bar\xi]_{\bar\partial}\mapsto[\bar\xi\wedge\alpha_t\wedge\beta_t\wedge\gamma_t\wedge\bar\gamma_t]_{\bar\partial}$ for $\xi\in\{\alpha_t,\,\beta_t\}$, while the isomorphism $ H^1_{DR}(X_t,\,\C)\simeq H^5_{DR}(X_t,\,\C)$ is canonically defined by $\{\zeta\}\mapsto\{\zeta\wedge\alpha_t\wedge\beta_t\wedge\gamma_t\wedge\bar\gamma_t\}$ for $ \zeta\in\{\bar\alpha_t,\,\bar\beta_t\}$ and by $\{\zeta\}\mapsto\{\zeta\wedge\gamma_t\wedge\bar\alpha_t\wedge\bar\beta_t\wedge\bar\gamma_t\}$ for $ \zeta\in\{\alpha_t,\,\beta_t\}$.  \hfill $\Box$

\vspace{3ex}

 Now, we recall two standard facts that prove between them that every elliptic curve (in particular, the fibre of the Albanese map $\pi:=\pi_0:X_0\to B_0$) is {\bf self-dual}.

\begin{Prop}\label{Prop:standard_elliptic-curves} (see e.g. [Dem97, $\S.10.2$]) Let $X$ be a compact complex manifold such that $\mbox{dim}_\C X=1$ (i.e. $X$ is a compact {\bf complex curve}). 

\vspace{1ex}

$(i)$\, The Jacobian torus $\mbox{Jac}(X)$ of $X$ coincides with its Albanese torus $\mbox{Alb}(X)$. Moreover, for every point $a\in X$, the Jacobi map

$$\Phi_a:X\longrightarrow \mbox{Jac}(X), \hspace{3ex} x\mapsto{\cal O}([x]-[a]),$$

\noindent coincides with the Albanese map

$$\alpha:X\longrightarrow \mbox{Alb}(X)=\mbox{Jac}(X).$$

$(ii)$\, If $X$ is an {\bf elliptic curve} (i.e. $g=1$, where $g:=h^{0,\,1}(X)$ is the genus of the complex curve $X$), then $\Phi_a=\alpha$ is an {\bf isomorphism}, i.e.

$$X\simeq \mbox{Jac}(X) = \mbox{Alb}(X).$$

\noindent In particular, since the dual tori $\mbox{Jac}(X)$ and $\mbox{Alb}(X)$ coincide, $X$ is self-dual.

\end{Prop}

\vspace{3ex}

 We can now infer the main result of this paper showing that the Iwasawa manifold is its own dual in a simple sense pertaining to its Albanese torus and map. This self-duality point of view complements those considered in [Pop17].

\begin{The}\label{The:Iwasawa_mirror_self-duality}  The Iwasawa manifold $X=X_0$ is {\bf its own dual} in the sense that in its Albanese map description 

$$\pi=\pi_0:X_0\longrightarrow B_0:=\mbox{Alb}(X_0)$$

\noindent as a locally holomorphically trivial fibration by elliptic curves $\C/\Z[i]$ over the $2$-dimensional complex torus $\C/\Z[i]\times\C/\Z[i]$, both the base $\mbox{Alb}(X_0)$ and the fibre $\pi_0^{-1}(s)$ are (sesquilinearly) self-dual tori.

\end{The}

\noindent {\it Proof.} The self-duality of $\mbox{Alb}(X_0)$ was proved in Lemma \ref{Lem:identification_tori}, while the self-duality of $\pi_0^{-1}(s)$ is the standard fact recalled in Proposition \ref{Prop:standard_elliptic-curves}.  \hfill $\Box$

\vspace{6ex}

\noindent {\bf References.} \\

\noindent [Ang14]\, D. Angella --- {\it Cohomological Aspects in Complex Non-K\"ahler Geometry} --- LNM 2095, Springer (2014).

\vspace{1ex}

\noindent [Dem 97]\, J.-P. Demailly --- {\it Complex Analytic and Algebraic Geometry}---http://www-fourier.ujf-grenoble.fr/~demailly/books.html

\vspace{1ex}

\noindent [Pop13]\, D. Popovici --- {\it Holomorphic Deformations of Balanced Calabi-Yau $\partial\bar\partial$-Manifolds}--- arXiv e-print AG 1304.0331v1.

\vspace{1ex}

\noindent [Pop14]\, D. Popovici --- {\it Deformation Openness and Closedness of Various Classes of Compact Complex Manifolds; Examples} --- Ann. Sc. Norm. Super. Pisa Cl. Sci. (5), Vol. XIII (2014), 255-305.

\vspace{1ex}

\noindent [Pop15]\, D. Popovici --- {\it Aeppli Cohomology Classes Associated with Gauduchon Metrics on Compact Complex Manifolds} --- Bull. Soc. Math. France {\bf 143} (3), (2015), p. 1-37.

\vspace{1ex}

\noindent [Pop17]\, D. Popovici --- {\it Non-K\"ahler Mirror Symmetry of the Iwasawa Manifold} --- arXiv e-print AG 1706.06449v1.

\vspace{1ex}

\noindent [PU14]\, D. Popovici, L. Ugarte --- {\it The sGG Class of Compact Complex Manifolds} --- arXiv e-print DG 1407.5070v1.

\vspace{1ex}

\noindent [Uen75]\, K. Ueno --- {\it Classification Theory of Algebraic Varieties and Compact Complex Spaces} --- LNM {\bf 439} (1975).

\vspace{1ex}

\noindent [Sch07]\, M. Schweitzer --- {\it Autour de la cohomologie de Bott-Chern} --- arXiv e-print math. AG/0709.3528v1.

\vspace{6ex}

\noindent Institut de Math\'ematiques de Toulouse, Universit\'e Paul Sabatier,

\noindent 118 route de Narbonne, 31062 Toulouse, France

\noindent Email: popovici@math.univ-toulouse.fr

\end{document}